\documentclass[11pt,fleqn,twoside]{article}
\usepackage{amsfonts,latexsym}
\makeatletter
\newcommand{\prava}{\footnotesize\it
\begin{flushright}
\begin{minipage}{18cm}
Copyright \copyright 1998 by M.L. Gandarias, P. Venero, and J. Ramirez
\end{minipage}
\end{flushright}}

\newcommand{\name}[1]{\begin{flushleft}
                       \LARGE \bf #1
                       \end{flushleft}\vspace{-3mm}}

\newcommand{\Author}[1]{\begin{flushleft}
                       \it #1 \end{flushleft}}

\newcommand{\Adress}[1]{\begin{flushleft}
                       \it #1 \end{flushleft}}

\newcommand{\Date}[1]{\begin{flushleft}
                      \small  \it #1 \end{flushleft}}

\newcommand{\ehkol}{Author \ name}
\newcommand{\ohkol}{Article \ name}
\renewcommand{\@evenhead}{
\hspace*{-3pt}\raisebox{-15pt}[\headheight][0pt]{\vbox{\hbox to \textwidth
{\thepage \hfil \ehkol}\vskip4pt \hrule}}}
\renewcommand{\@oddhead}{
\hspace*{-3pt}\raisebox{-15pt}[\headheight][0pt]{\vbox{\hbox to \textwidth
{\ohkol \hfil \thepage}\vskip4pt\hrule}}}
\renewcommand{\@evenfoot}{}
\renewcommand{\@oddfoot}{}

\newcommand{\be}{\begin{equation}}
\newcommand{\ee}{\end{equation}}
\newcommand{\ba}{\hspace*{-5pt}\begin{array}}
\newcommand{\ea}{\end{array}}

\newcommand{\ds}{\displaystyle}
\makeatother

\begin{document}
\setcounter{page}{234}
\thispagestyle{empty}

\renewcommand{\ehkol}{M.L. Gandarias, P. Venero and J. Ramirez}
\renewcommand{\ohkol}{Similarity Reductions for a Nonlinear Dif\/fusion
Equation}

\begin{flushleft}
\footnotesize {\sf
Journal of Nonlinear Mathematical Physics \qquad 1998, V.5, N~3},\
\pageref{gandarias-fp}--\pageref{gandarias-lp}.
\hfill
{{\sc Letter}}
\end{flushleft}

\vspace{-5mm}

\renewcommand{\footnoterule}{}
{\renewcommand{\thefootnote}{}
 \footnote{\prava}}

\name{Similarity Reductions for \\ a Nonlinear Dif\/fusion
Equation}\label{gandarias-fp}

\Author{M.L. GANDARIAS, P. VENERO and J. RAMIREZ}

\Adress{Departamento de Matematicas, Universidad de Cadiz, \\
P.O. Box 40, 1510 Puerto Real, Cadiz, Spain\\
E-mail: mlgand@merlin.uca.es}

\Date{Received November 19, 1997, Revised May 27, 1998, Accepted July
8, 1998}

\begin{abstract}
\noindent
Similarity reductions and new exact solutions are obtained for
a nonlinear dif\/fusion equation. These are obtained by using the classical
symmetry group and reducing the partial dif\/ferential equation to
various ordinary dif\/ferential equations. For the equations so
obtained, f\/irst integrals are deduced which consequently give rise to
explicit solutions. Potential symmetries, which are realized as local
symmetries of a related auxiliary system, are obtained. For some
special nonlinearities new symmetry reductions and exact solutions
are derived by using the nonclassical method. 
\end{abstract}

\section{Introduction}
Equation
\begin{equation}\label{gandarias:eq6}
u_t =\left(u^ n\right)_{xx}+\frac{C}{x+\lambda}\left(u^ n\right)_x,
\end{equation}
with
$n \in {\mathbb Q} \backslash \{0,1\}$ and $\lambda \in {\mathbb R}$,
corresponds to nonlinear dif\/fusion with convection. When $C=0$
equation (\ref{gandarias:eq6}) becomes the one-dimensional porous medium
equation
\begin{equation}\label{gandarias:eq4}
u_t =\left(u^ n\right)_{xx}.
\end{equation}
 A complete group classif\/ication for  (\ref{gandarias:eq4}) was derived by Ovsiannikov
 \cite{gandarias:ovsi}  and by  Bluman \cite{gandarias:blu0,gandarias:blucoll}. A classif\/ication for
 Lie-B\"acklund symmetries was obtained by Bluman and Kumei
  \cite{gandarias:bluku0}.
 The main known exact solutions of nonlinear dif\/fusion (\ref{gandarias:eq4}) are
 summarized by Hill \cite{gandarias:hi}. In
 \cite{gandarias:hi,gandarias:hi1,gandarias:hi2}, Hill {\em et al}
  deduced a number of f\/irst integrals for stretching similarity solutions
  of the nonlinear dif\/fusion equation, and of general high-order nonlinear
  evolution equations, by two dif\/ferent  integration procedures.
King \cite{gandarias:ki1} obtained approximate solutions to the porous medium equation
(\ref{gandarias:eq4}).

The basic idea of any similarity solution is to assume a functional
form of the solution which enables a PDE to be reduced to an ODE. The
majority of known exact solutions of (\ref{gandarias:eq4}) turn out to be
similarity solutions, even though originally they might have been
derived, say by a separation of variable technique or as traveling
wave solutions. For 
\[
u_t=\left(u^n\right)_{xx}+C\left(u^n\right)_x,
\]
which is the Boussinesq equation of hydrology involved in various
f\/ields of petroleum technology and ground water hydrology, several exact
solutions have been obtained by using isovector method \cite{gandarias:bhuta}.

More often than not the spatial dependent factors are assumed to be
constant, although there is no fundamental reason to assume so.
Actually, allowing for their spatial dependence enables one to
incorporate additional factors into the study which may play an
important role. For instance, in a porous medium this may account for
intrinsic factors, like medium contamination with another material,
or in plasma, this may express the impact that solid impurities
arising from the walls have on the enhancement of the radiation
channel. Knowing the importance of the ef\/fect of space-dependent
parts on the overall dynamics of the nonlinear dif\/fusion equation, a
group classif\/ication for
\begin{equation}\label{gandarias:poro}
u_t=\left(u^n\right)_{xx}+ f(x)u^su_x+g(x)u^m
\end{equation}
was derived in \cite{gandarias:ml}, by studying those spatial forms which admit
the classical symmetry group. Both the symmetry group and the spatial
dependence was found through consistent application of the Lie-group
formalism. The behaviour of the interface of a related problem with
(\ref{gandarias:eq6}) has been studied by Okrasinski \cite{gandarias:okra}.

In this work we use the invariance of equation (\ref{gandarias:eq6}) under
one-parameter group of transformations to reduce the PDE (\ref{gandarias:eq6})
to various ordinary dif\/ferential equations. Most of the required
theory and description of the method can be found in
\cite{gandarias:bluku1,gandarias:hill,gandarias:olv,gandarias:ovsi,gandarias:steph}.
Following Hill \cite{gandarias:hi,gandarias:hi1,gandarias:hi2},
we have deduced some exact solutions of equation (\ref{gandarias:eq6}) by fully
integrating the ODE's derived. 

An obvious limitation of group-theoretic methods based in local
symmetries, in their utility for particular PDE's, is that many of
these equations do not have local symmetries. It turns out that PDE's
can admit nonlocal symmetries whose inf\/initesimal generators depend
on integrals of the dependent variables in some specif\/ic manner.

In \cite{gandarias:bluku0,gandarias:bluku1} Bluman introduced a method to f\/ind a new
class of symmetries for a PDE. By writing a given PDE, denoted by
R\{x,t,u\}, in a conserved form, a 
related system denoted by S\{x,t,u,v\} with potentials as additional
dependent variables, is obtained. Any Lie group of point
transformations admitted by S\{x,t,u,v\} induces a symmetry for
R\{x,t,u\}; when at least one of the generators of the group depends
explicitly on the potential; then the corresponding symmetry is
neither a point nor a Lie-B\"acklund symmetry. These symmetries of
R\{x,t,u\} are called {\em potential} symmetries.

The nature of potential symmetries allows one to extend the uses of
point symmetries to such nonlocal symmetries. In particular:
Invariant solutions of S\{x,t,u,v\} yield solutions of R\{x,t,u\}
which are not invariant solutions for any local symmetry admitted by
R\{x,t,u\}. Potential symmetries for equation (\ref{gandarias:poro}), when it
can be written in a conserved form, have been recently derived
\cite{gandarias:ml1}.

In order to f\/ind potential symmetries of (\ref{gandarias:eq6}), we write this
equation in the conserved form 
\begin{equation}\label{gandarias:conserv}
D_xF-D_tG.
\end{equation}
The associated auxiliary system S\{x,t,u,v\} is then given by
\begin{equation}\label{gandarias:sis}
v_x=(x+\lambda)u, \qquad v_t=(x+\lambda)\left(u^n\right)_x+(C-1) u^n.
\end{equation}
 Suppose S\{x,t,u,v\} admits a local Lie group of
transformations with inf\/initesimal generator
\begin{equation}\label{gandarias:ge1}
X_S=p(x,t,u,v)\frac{\partial}{\partial
x}+q(x,t,u,v)\frac{\partial}{\partial t}+r(x,t,u,v)
\frac{\partial}{\partial u}+s(x,t,u,v)\frac{\partial}{\partial v}.
\end{equation} 
This group maps any solution of S\{x,t,u,v\} to another solution of
S\{x,t,u,v\} and hence induces a mapping of any solution of
R\{x,t,u\} to another solution of R\{x,t,u\}. 
Thus (\ref{gandarias:ge1}) def\/ines a symmetry group of R\{x,t,u\}. If
\begin{equation}\label{gandarias:cod}
\left(\frac{\partial p}{\partial x}\right)^2
+\left(\frac{\partial q}{\partial v}\right)^2
+\left(\frac{\partial r}{\partial v}\right)^2\neq
0
\end{equation}
then (\ref{gandarias:ge1}) yields a nonlocal symmetry of R\{x,t,u\}. Such
nonlocal symmetry is called a {\em potential} symmetry of R\{x,t,u\}
\cite{gandarias:bluku0,gandarias:bluku1}.

Motivated by the fact that symmetry reductions for many PDE's are
known that are not obtained by using the classical Lie group method,
there have been several generalizations of the classical Lie group
method for symmetry reductions. Bluman and Cole [3] introduced the
nonclassical method to study the symmetry reductions of the heat
equation and Clarkson and Mansf\/ield \cite{gandarias:cla3} presented an algorithm for
calculating the determining equations associated with the
nonclassical method. The basic idea of the method is that the PDE
(\ref{gandarias:eq6}) is augmented with the invariance surface condition
\begin{equation}\label{gandarias:sur}
pu_x+qu_t-r=0,
\end{equation}
which is associated with the vector f\/ield
\begin{equation}\label{gandarias:vect}
V=p(x,t,u) \frac{\partial}{\partial x}
+q(x,t,u)\frac{\partial}{\partial t}+r(x,t,u)
\frac{\partial}{\partial u}.
\end{equation} 
By requiring that both (\ref{gandarias:eq6}) and (\ref{gandarias:sur})
are invariant under the transformation with in\-f\/i\-ni\-te\-si\-mal generator
(\ref{gandarias:vect}) one obtains an
overdetermined nonlinear system
of equations for the in\-f\/i\-ni\-te\-si\-mals $p(x,t,u),$ $q(x,t,u),$ $r(x,t,u).$
The number of determining equations arising in the nonclassical
method is smaller than for the classical method, consequently the
set of solutions is in general, larger than for the classical
method as in this method one requires only the subset of
solutions of (\ref{gandarias:eq6}) and (\ref{gandarias:sur})
to be invariant under the inf\/initesimal generator (\ref{gandarias:vect}).
However, the associated vector f\/ields do not form a vector space.
These methods were generalized and called conditional symmetries by
Fushchych and Nikitin \cite{gandarias:fu} and
also by Olver and Rosenau \cite{gandarias:olv2,gandarias:olv3} to include ``weak symmetries'',
``side conditions'' or ``dif\/ferential constraints''.

The nonclassical symmetries of the nonlinear dif\/fusion equation
(\ref{gandarias:poro}) with an absorption term and $f(x)=0$, as well as new
exact solutions were derived in \cite{gandarias:ml2}. In this work we obtain
the special values of the parameter $n$ such that
nonclassical symmetries for (\ref{gandarias:eq6}) can be derived. We also report the
reduction obtained as well as some new exact solutions.

\section{Classical symmetries: Exact solutions}
In this case, we f\/ind that the most general Lie group of point
transformations admitted by (\ref{gandarias:eq6}) is:
\begin{enumerate} 
\item
For $n\neq 0,1$ and $\ds C\neq \frac{3n+1}{n+1}$ we obtain a
three-parameter group ${\cal G}_3.$ Associated with this Lie group is
its Lie algebras ${\cal L}_3,$ which can be respectively represented
by the set of all the generators $\{V_i\}_{i=1}^3.$ These 
generators are:
\begin{equation}\label{gandarias:ge2}
\frac{\partial}{\partial t}, \qquad
V_2=(x+\lambda)\frac{\partial}{\partial x}+2t \frac{\partial}{\partial t},
\qquad
 V_3=-t \frac{\partial}{\partial t}
+\frac{u}{n-1}
\frac{\partial}{\partial u}.
\end{equation}
\item
For $n\neq -1,0,1$ and $\ds C=\frac{3n+1}{n+1}$, we obtain a
four-parameter group ${\cal G}_4$. Its inf\/initesimal generators
$\{V_i\}_{i=1}^4$ are: 
\begin{equation}\label{gandarias:ge22}
\ba{l}
V'_1=V_1, \qquad V'_2=V_2, \qquad V'_3=V_3, 
\vspace{3mm}\\
\ds V'_4= \frac{(x+\lambda)^{2-C}}{(C-2)(C-1)}
\frac{\partial}{\partial x}-
\frac{2(x+\lambda)^{1-C}}{(C-1)(n-1)}u
\frac{\partial}{\partial u}.
\ea
\end{equation}
\end{enumerate}
To ensure that an optimal set of reductions is obtained from the
symmetries of (\ref{gandarias:eq6}), the optimal system is determined for
$\lambda\neq 0.$ The case $\lambda =0$ was derived in \cite{gandarias:ml}. In
Table 1, we list the nontrivial optimal system $\{U_i\}$ with
$i=1,2,3,4.$ We also list the corresponding similarity variables and 
similarity solutions.

\medskip

\begin{center}
\begin{tabular}{|llll|}
\multicolumn{4}{l}{\parbox{13cm}{
{\bf Table 1: }{\small Each row show the inf\/initesimal
generators of the optimal system,
the corresponding similarity variables and similarity solutions.}}}\\[4mm]
\hline
  $i$ & $U_i$ & $z_i$ & $u_i$ \\ \hline
 \mbox{ }  & \mbox{ } & \mbox{ } & \mbox{ } \\
 $1$ & $V'_2+aV'_3$ & $\ds \frac{(x+\lambda)^{2-a}}{t}$ &
$ h(z)(x+ \lambda)^{{{a}\over{n-1}}}$ \\
\mbox{ } & \mbox{ } & \mbox{ } & \mbox{ } \\
 $2$ & $aV'_1+V'_2+2V'_3$ & $ e^ {-at}(x+ \lambda)$ &
$ h(z)(x+\lambda)^{{{2}\over{n-1}}}$  \\
\mbox{ } & \mbox{ }& \mbox{ } & \mbox{ } \\
$3$ & $aV'_3+V'_4$ & $ {\rm c_1}
\left(n-1\right)\left(x+\lambda\right)^{C-1}+
\log t $ & $\ds  {{h(z)e^{{\rm c_1}
\left(x+\lambda\right)^{{{2n}\over{n+1}}}}}
\over{\left(x+\lambda\right)^{{{2}\over{n+1}}}}} $\\
\mbox{ } & \mbox{ } & \mbox{ } & \mbox{ } \\
$4$ & $aV'_1+V'_4$ & $\ds  \frac{a(C-2)(x+\lambda)^{C-1}}{t}$ &
$h(z)[(n-1)(x+\lambda)]^{\frac{2(2-C)}{n-1}}$\\
 \mbox{ } & \mbox{ } & \mbox{ } & \mbox{ } \\ \hline
\end{tabular} 
\end{center}

\medskip

In Table 1, case $i=1$, the constant $a\in { \mathbb Q} \backslash
\left\{2\right\}$
and, in case $i=3$, $\ds c_1=\frac{a}{n+1}$.

The ODE to which the PDE (\ref{gandarias:eq6}) is reduced by means of the inf\/initesimal
generator $U_1$ is
\[
\ba{l}
\ds
 {{d h}\over{d z}}\left({{1}\over{\left(a-2\right)^{2}h^{n-1}
 n}}-{{C_{1}+n+a-1}\over{\left(a-2\right)\left(n-1\right)z
 }}\right)
 \vspace{3mm}\\
 \ds \qquad +{{ah\left(C_{1}-n+1\right)}\over{\left(a-2
 \right)^{2}\left(n-1\right)^{2}z^{2}}}+
 \frac{n-1}{h}\left(\frac{dh}{dz}\right)^2+{{d^{2} h}\over{d z^{2}}}=0,
\ea
\]
with
\[
C_1=Cn+an-C,
\]
after taking $\ds h=y^{\frac{1}{n}}$ it becomes
\[
 -{{d y}\over{d z}}\left({{C_1+n+a-1}\over{\left(a-2
 \right)\left(n-1\right)z}}-{{y^{{{1}\over{n}}-1}}\over{\left(a-2
 \right)^{2}n}}\right)+{{an\left(C_1-n+1\right)y
 }\over{\left(a-2\right)^{2}\left(n-1\right)^{2}z^{2}}}+{{d^{2} y
 }\over{d z^{2}}} =0.
\]
In particular, the second order nonlinear dif\/ferential equation for
 $y(z)$ obtained for $i=1$,  admits  f\/irst integrals for  some values
 of $a$. As an example, we consider
 \[
 a=2-2n.
 \]
 The f\/irst integral  for $a=2-2n$ is
\[
 -{{y}\over{z}}+{{Cy}\over{2nz}}-
{{y}\over{2nz}}
+{{{{d y}\over{d z}}}}+
{{y^{{{1}\over{n}}}}\over{4n^{2}}}=c_1,
\]
where $c_1$ is the integration constant.
If this constant is zero  we obtain  for $a=2-2n$
\[
y= {{\left(\left(1-n\right)z^{{{4n^{2}+{\rm k_2}}\over{4n^{2}}}}
 +\left(4{\rm k_1}n^{3}+{\rm k_1}{\rm k_2}n-{\rm k_1}{\rm k_2}
 \right)z^{{{{\rm k_2}}\over{4n^{3}}}}\right)^{{{n}\over{n-1}}}
 }\over{\left(4n^{3}+{\rm k_2}n-{\rm k_2}\right)^{{{n}\over{n-1}}}
z^{{{{\rm k_2}}\over{4n^{2}-4n}}}}},
\]
with
\[
 k_2=\left(2c-2\right)n-4n^{2},
\]
\[
h= {{\left(\left(1-n\right)z^{{{4n^{2}+{\rm k_2}}\over{4n^{2}}}}
 +\left(4{\rm k_1}n^{3}+{\rm k_1}{\rm k_2}n-{\rm k_1}{\rm k_2}
 \right)z^{{{{\rm k_2}}\over{4n^{3}}}}\right)^{{{1}\over{n-1}}}
 }\over{\left(4n^{3}+{\rm k_2}n-{\rm k_2}\right)^{{{1}\over{n-1}}}
z^{{{{\rm k_2}}\over{4n^{3}-4n^{2}}}}}}.
\]
 Substituting in the similarity solution we obtain the exact solution
\[
u= {{\left(\left(n-1\right)z^{{{4n^{2}+{\rm k_2}}\over{4n^{2}}}}
 +\left(-4{\rm k_1}n^{3}-{\rm k_1}{\rm k_2}n+{\rm k_1}{\rm k_2}
 \right)z^{{{{\rm k_2}}\over{4n^{3}}}}\right)^{{{1}\over{n-1}}}
 }\over{\left(-4n^{3}-{\rm k_2}n+{\rm k_2}\right)^{{{1}\over{n-1}}}
z^{{{{\rm k_2}}\over{4n^{3}-4n^{2}}}}}}.
\]
When $C=0$ this is the well-known dipole solution \cite{gandarias:hi}.

In particular, the second order nonlinear dif\/ferential equation for
$h(z)$ obtained for $i=2,$ is
\[
 {{d h}\over{dz}}\left({{Cn+
4n-C}\over{\left(n-1\right)z}}
 +{{a}\over{h^{n-1}nz}}\right)+{{2h\left(Cn+n-C+1\right)
 }\over{\left(n-1\right)^{2}z^{2}}}+
 \frac{n-1}{h}\left(\frac{dh}{dz}\right)^2
 +{{d^{2} h}\over{d z^{2}}} =0.
 \]
For the special value
\[
C= -{{n+1}\over{n-1}},
\]
we integrated once to obtain
\[
 h^{n-1}{{d h}\over{d z}}z+{{h^{n}\left(3n-1\right)}\over{
 \left(n-1\right)n}}-{{h^{n}}\over{n}}+{{ah}\over{n}}=c_1.
\]
After taking
\[
h= y^{{{1}\over{n-1}}}
\]
it becomes
\[
{{2y}\over{z}}-{{a}\over{nz}}+{{a}\over{z}}+{{d y}\over{d z}} =c_1,
\]
where $c_1$ is the integration constant. If this constant is zero  we
obtain
\[
y= {{c}\over{z^{2}}}-{{a\left(n-1\right)}\over{2n}},
\]
so that
\[
h= \left({{ce^{2at}}\over{x^{2}}}-{{a\left(n-1\right)}\over{2
 n}}\right)^{{{1}\over{n-1}}}x^{{{2}\over{n-1}}}.
\]
Substituting in the similarity solution we obtain the exact solution
\[
u= \left({{ce^{2at}}\over{x^{2}}}-{{a\left(n-1\right)}\over{2
 n}}\right)^{{{1}\over{n-1}}}x^{{{2}\over{n-1}}}.
\]

The second order nonlinear ODE obtained, after taking
$h(z)=y^{1/n}$, for $i=3$ is
\[
-{{\left(n+1\right)^{2}y^{{{1}\over{n}}-1}e^ {-{ z }}}\over{4{\rm c_1}^{2}
\left(n-1\right)^{2}n^{3}}}{{d y}\over{d z}}+{{d ^{2} y}\over{d z^{2}}}
+{{2n}\over{n-1}}{{d y}\over{d z}}+{{n^{2}
 y}\over{\left(n-1\right)^{2}}} =0,
\]
where $\ds c_1=\frac{a}{n+1}$.

The second order nonlinear ODE obtained, after taking
$h(z)=y^{1/n}$, for $i=4$, is
\[
\frac{d^2y}{dz^2}
+\left(\frac{1}{2z}+c_1y^{1/n-1}\right)\frac{dy}{dz}=0,
\]
where $\ds c_1=\frac{(n-1)^{\frac{n-3}{n+1}}(n+1)^3}{16an^3}$.

\section{Potential Symmetries}

In order to f\/ind potential symmetries of (\ref{gandarias:eq6}),
we write this equation in the conserved form~(\ref{gandarias:conserv}), where
\[
G=(x+\lambda)u,
\]
\[F=(x+\lambda)\left(u^n\right)_x+(C-1)u^n.
\]
The associated auxiliary system is given by (\ref{gandarias:sis}). Besides
$\ds X_1=\frac{\partial}{\partial t}$ and
$\ds X_2=\frac{\partial}{\partial v}$ we obtain $X_3,\ldots,X_6$,
given in Table 2.

\medskip

\begin{center}
\begin{tabular}{|lllllll|}
\multicolumn{7}{l}{\parbox{13cm}{
{\bf Table 2:}
{\small
Each row show the inf\/initesimal generators,
the corresponding similarity variables and similarity solutions.}}}\\[3mm]
 \hline
$(a)$ & $p$ & $q$ & $r$ & $s$ & $n$ & $C$\\ \hline
\mbox{ } & \mbox{ }& \mbox{ }& \mbox{ }& \mbox{ }& \mbox{ } & \mbox{ }\\
$X_3$ & $ {x+\lambda}$ & $0$ & $\ds \frac{2u}{n-1}$ &
$\ds \frac{2nv}{n-1}$ & $\neq 0,1$ &
{\small arbitrary} \\
\mbox{ } & \mbox{ }& \mbox{ }& \mbox{ }& \mbox{ }& \mbox{ }& \mbox{ }\\
$X_4$ & $(x+\lambda)$ & $2nt$ & $-2u$ & $0$ & ${\neq 0, 1}$ & {\small
 arbitrary} \\
\mbox{ } & \mbox{ }& \mbox{ }& \mbox{ }& \mbox{ }& \mbox{ } & \mbox{ }\\
$X_5$ & $\ds k(x+\lambda)^{\frac{1-n}{1+n}}$ & $0$ &
$\ds -\frac{2k}{n+1}(x+\lambda)^{-\frac{2n}{n+1}}$ & $0$ & $\neq -1,0,1$ &
$\ds \frac{3n+1}{n+1}$\\
\mbox{ } & \mbox{ }& \mbox{ }& \mbox{ }& \mbox{ }& \mbox{ } & \mbox{ }\\
$X_6$ & $2(x+\lambda)v$ & $0$ & $-2(x+\lambda)^2u^2-2uv$ & $v^2$ & $-1$ &
$\ds -\frac{1}{2}$\\
\mbox{ } & \mbox{ }& \mbox{ }& \mbox{ }& \mbox{ }& \mbox{ } & \mbox{ }\\
 \hline
\end{tabular}
\end{center}

\medskip

$X_3,$ $X_4$ and $X_5$ project onto point
symmetries of (\ref{gandarias:eq6});
while $X_6$ induces potential symmetries admitted by
(\ref{gandarias:eq6}).

Solving the characteristic equation,
we obtain the similarity variable $z=t$ and
similarity solution $v=\sqrt{x+\lambda}E(t)$.
In this case, $E(t)$ satisf\/ies
$E'(t)=0$, so $E={\rm const}$ and we obtain the trivial solution
\[
v=C_1\sqrt{x+\lambda}, \qquad u=\frac{C_1}{2(x+\lambda)^\frac{3}{2}}.
\]

We must note that, although the inf\/initesimal
$p$ depends explicitly on $v$, the
similarity variable does not depend on $v$.

\section{Nonclassical symmetries}
To apply the nonclassical method to (\ref{gandarias:eq6})
we require  (\ref{gandarias:eq6}) and (\ref{gandarias:sur}) to be invariant
under the inf\/initesimal generator (\ref{gandarias:vect}).
In the case $q \neq 0,$ without loss of generality,
we may set $q(x,t,u)=1$. The nonclassical method applied to
(\ref{gandarias:eq6}) give rise to four determining equations for the inf\/initesimals.
\[
{{\partial^{2} p}\over{\partial u^{2}}}u-n{{\partial p}\over{\partial u}}
+{{\partial p}\over{\partial u  }} =0,
\]
\[
\ba{l}
\ds -\left(n{{\partial^{2} r}\over{\partial u^{2}}}-2n{{\partial^{2} p}
\over{\partial u \partial x  }}+2\frac{nC}{x+\lambda}{{\partial p}
\over{\partial u}}\right)u^{n-1}
\vspace{3mm}\\
\ds \qquad -\left(n-1\right)n{{\partial  r}\over{\partial u}}
u^{n-2}+\left(n-1\right)nru^{n-3}-2p{{\partial p
 }\over{\partial u}}=0,
 \ea
\]
\[
\ba{l}
\ds -\left(2n{{\partial^{2} r}\over{\partial u
\partial x}}-n{{\partial^{2} p}\over{\partial x^{2} }}
+\frac{nC}{x+\lambda}{{\partial p}\over{\partial x}}+
\frac{-nC}{(x+\lambda)^2}p\right)u^{n+1}
\vspace{3mm}\\
\ds \qquad -2\left(n-1\right)n{{\partial r}\over{\partial x}}u^{n}+
\left(2{{\partial p}\over{\partial  u}}r-2p{{\partial p}
\over{\partial x}}-{{\partial p}\over{\partial t}}\right)u^{2}+
 \left(n-1\right)pru=0,
 \ea
 \]
\[
 -\left(n{{\partial^{2} r}\over{\partial x^{2}}}+
\frac{nC}{x+\lambda}
{{\partial r}\over{\partial x}}\right)u^{n}+
\left({{\partial r}\over{\partial t}}+
2{{\partial p}\over{\partial x}}r\right)u-
 \left(n-1\right)r^{2}=0,
\]
where $\ds f(x)=\frac{nC}{x+\lambda}.$
Solutions of this system depend in a fundamental way on the values of~$n$.
By solving the determining equations  we obtain
\[
p= {p_2(x,t)}u^{n}+{p_1(x,t)}.
\]
We can distinguish now the following: if $\ds n\in
\left\{0,-1,-\frac{1}{2}\right\}$
 we recover the classical symmetries, and if
 $\ds n \notin \left\{ 0,-1,-\frac{1}{2}\right\},$  we obtain that
\[
 r={\rm a_2}u^{n+2}+{\rm a_3}u^{n+1}+{\rm a_1}u^{n+1}+{{{\rm r_2}
 }\over{u^{n-1}}}+{\rm a_4}u^{2}+{\rm r_1}u,
 \]
 where
\[
 a_1= -\frac{Cp_2}{n(x+\lambda)}, \qquad
a_2= -{2{\rm p_2}^{2}
\over{\left(n+1\right)\left(2n+1\right)}},
\qquad a_3= \frac{1}{n}\frac{\partial p_2}{\partial x} ,
\qquad a_4= -{2{\rm p_1}{\rm p_2}\over{n+1}}
\]
and $p_1,$ $p_2,$ $r_1,$ and $r_2$  are related by two  conditions.
After considering the special values for $n$ for which
new symmetries dif\/ferent from Lie classical symmetries
can be obtained,  we can now distinguish the following:
\begin{itemize}
\item
for $\ds n\neq \frac{1}{2}$ we recover the classical symmetries,
\item
for $\ds n=\frac{1}{2}$
it follows that $p_2=0$ and $p_1,$ $r_1,$ and $r_2$
are related by the following conditions
\[-{\rm p_1}{\rm r_1}-4{\rm p_1}
{{\partial {\rm p_1}}\over{\partial x}}-2
{{\partial {\rm p_1}}\over{\partial t}}=0,
\]
\[
 -{\rm p_1}{\rm r_2}-{{\partial {\rm r_1}}\over{\partial x}}+
 {{\partial^{2} {\rm p_1}}\over{\partial x^{2}}}
-2\frac{nC}{x+\lambda}\frac{\partial p_1}{\partial x}
+2\frac{nC}{(x+\lambda)^2}p_1=0,
\]
\[
2{{\partial {\rm r_1}}\over{\partial t}}+
{\rm r_1}^{2}+4{{\partial {\rm p_1}}\over{\partial x}}{\rm r_1}=0,
\]
\[
 2{{\partial {\rm r_2}}\over{\partial t}}+2{\rm r_1}{\rm r_2}+
 4{{\partial {\rm p_1} }\over{\partial x}}{\rm r_2}-
 {{\partial^{2} {\rm r_1}}\over{\partial x^{2}}}-
 2\frac{nC}{x+\lambda}{{\partial  {\rm r_1}}\over{\partial x}}=0,
\]
\[
 -{{\partial^{2} {\rm r_2}}\over{\partial x^{2}}}-
 2\frac{nC}{x+\lambda}{{\partial{\rm r_2}}\over{\partial x}}+
 {\rm r_2}^{2}=0.
\]
We do not solve the above equations in general, but consider some
special solutions:
\end{itemize}
\begin{enumerate}
\item
Choosing $p_1=k,$  $r_1=0$ and  $C=2$, then
$\ds r_2= \frac{2}{(x+\lambda)^2}$
and we obtain the nonclassical ansatz
\[
z=x-kt, \qquad  u=\left(h(z)-\frac{1}{k(x+\lambda)}
\right)^2,
\]
where  $k \neq 0$ and  $h(z)$ satisf\/ies the following ODE
\begin{equation}\label{gandarias:od}
h''+2khh'=0;
\end{equation}
whose solutions are
\[
\ba{lllll}
\ds h(z)= \frac{k_4\tanh(k_4(z+k_2))}{k},
 & \mbox{if } & kk_3>0 & \mbox{and } & k_4=\sqrt{kk_3},
 \vspace{3mm}\\
\ds   h(z)=-\frac{k_4\tan (k_4(z+k_2))}{k},
& \mbox{if } & kk_3<0 & \mbox{and } &
k_4=\sqrt{-kk_3}.
\ea
\]

This leads to the exact solutions
\[
\ba{lll}
\ds u(x,t)=
 \left(\frac{k_4 \tanh (k_4(z+k_2))}
{k}-\frac{1}{k(x+C)}\right)^2 &
\mbox{if } & kk_3>0,
\vspace{3mm}\\
\ds u(x,t)=\left( \frac{ k_4 \tan
(k_4(z+k_2))}{k}-\frac{1}{k(x+C)}\right)^2, &
\mbox{if } & kk_3<0.
\ea
\]
\item
Choosing $p_1=p_1(x),$ $r_1=r_1(x)$ and $r_2=r_2(x)$ we obtain
\[
r(x,t,u)=-4p'u+\left[\frac{5p''_1}{p_1}-
\frac{Cp'_1}{(x+\lambda)p_1}+
\frac{C}{(x+\lambda)^2}\right]\sqrt{u},
\]
where $p_1$ must satisfy the following equations
\[
\frac{C}{2(x+\lambda)}=\frac{3p'_1}{2p_1}
-\frac{p''_1}{2p'_1}-\frac{c_1}{2p_1p'_1},
\]
\[
\ba{l}
\ds
-\frac{2C}{x+\lambda}\left[2p_1^2p'''_1+4p_1p'_1p''_1-(p'_1)^3\right]
+\frac{cp_1}{(x+\lambda)^2}
\left[Cp_1p''_1+8p_1p''_1+2(p'_1)^2\right]
\vspace{3mm}\\
\ds \qquad -\frac{C}{(x+\lambda)^3}(3C-2)p_1^2p'_1
+\frac{3}{(x+\lambda)^4}(C-2)Cp_1^3-5p_1^2p''''_1
\vspace{3mm}\\
\qquad +10p_1p'_1p'''_1+30p_1(p''_1)^2-10(p'_1)^2p''_1=0.
\ea
\]
We observe that these condition are satisf\/ied for
$p_1=k_1(x+\lambda)^{\frac{C-1}{2}}$
if $\ds C\in\left\{\frac{5}{3},-1,3\right\}$.
\begin{itemize}
\item
 For $C=3,$ $p_1=k_1(x+\lambda)$ we recover classical
symmetries.
\item
 For $\ds C=\frac{5}{3},$ we obtain
 \[
 p_1=k_1(x+\lambda)^{1/3}, \qquad
  r=-\frac{4k_1u}{3(x+a)^{2/3}}.
\]
Hence a nonclassical ansatz is
\[
z=\frac{3k_1}{2}(x+\lambda)^{2/3}-kt,
\qquad u=(x+\lambda)^{-4/3}h(z),
\]
where $h(z)$ satisf\/ies the following ODE
\[
-2hh''+(h')^2+-4kk_1^2h^{3/2}h'=0,
\]
whose solutions are
\[\!\!
\ba{lllll}
\ds h(z)=\frac{k_2}{4k_1^2k}\tan ({k_1k_4(z+k_3)}) & \!
\mbox{if } kk_2>0, & \!k_4=\sqrt{kk_2} & \!\mbox{and} & k_1<0
\vspace{3mm}\\
\ds h(z)=\frac{k_2}{4k_1^2k}\tanh ({k_1k_4(z+k_3)}) &
\!\mbox{if } kk_2<0, & \!k_4=\sqrt{-kk_2} & \!\mbox{and} & k_1>0
\ea
\]
\item
For $C=-1,$ we obtain
\[
p_1=\frac{k_1}{x+\lambda},
\qquad r=\frac{4(k_1u+2\sqrt{u})}{(x+\lambda)^{2}}.
\]
Hence a nonclassical reduction is
\[
z=\frac{x^2+2\lambda x}{k_1}-t,
\qquad u=\frac{1}{k_1^2}(h(z)(x+\lambda)^2-2)^2,
\]
where $h(z)$ satisf\/ies (\ref{gandarias:od}) with $k=k_1$.
\end{itemize}
\item
Choosing $p_2=r_2=0$ and $p_1=p_1(x),$ we obtain $r_1=-4p'_1$ and
$\ds f=\frac{5p'_1-2k_1}{2p_1},$ where $p_1$ satisf\/ies
\[
p_1^2p''_1+2p_1(p'_1)^2-2k_1p_1p'_1-k_2p_1=0,
\]
which is a classical reduction \cite{gandarias:ml}.
\end{enumerate}

\section{Concluding Remarks}
In this paper we have used the invariance of (\ref{gandarias:eq6})
under  group of transformations to reduce~(\ref{gandarias:eq6}) to ODEs.
We desired to minimize the search for group-invariant
solutions to that of f\/inding non-equivalent branches of solutions,
consequently we have constructed all the invariant
solutions with respect to  the one-dimensional optimal
system of subalgebras, as well as all the ODEs to
which (\ref{gandarias:eq6}) is reduced. For the equations so obtained,
f\/irst integrals have been deduced which give rise to
explicit solutions.
Potential symmetries  as well as nonclassical symmetries were
used to obtain new solutions of (\ref{gandarias:eq6}). The new solutions
are unobtainable by Lie classical symmetries.

\section{Acknowledgments}
It is a pleasure  to thank Professor Philip Rosenau for bringing the
porous medium equation to my attention,
as well as for his support on this work.
I am also grateful to Professor Romero for his helpful comments.

\label{gandarias-lp}

\begin{thebibliography}{44}
\footnotesize


\bibitem{gandarias:bhuta} Bhutani O.P.  and Vijayakumar K.,
  {Int. J. Engng. Sci.,} 1990, V.28, 375--387.
\bibitem{gandarias:blu0} Bluman G.W., Construction of solutions to partial
dif\/ferential equations by the use of transformations groups,
 Ph.D. Thesis, Inst. of Technology, California, 1967.
\bibitem{gandarias:bluu} Bluman G.W. and Cole J.D.,  {\it J. Math. Mech.,}
1969, V.18, 1025.
\bibitem{gandarias:blucoll} Bluman G.W. and  Cole J.D., Similarity Methods for
Dif\/ferential Equations, Springer, Berlin, 1974.
\bibitem{gandarias:bluku0} Bluman G.W. and Kumei S., {\it J. Math. Phys.,}
1980, V.21, 1019.
\bibitem{gandarias:bluku1} Bluman G.W. and  Kumei S.,
Symmetries and Dif\/ferential Equations, Springer, Berlin, 1989.

\bibitem{gandarias:cla} Clarkson P.A. and Kruskal M.D., {\it J.Math. Phys.,}
1989, V.30, 2201.

\bibitem{gandarias:cla3}  Clarkson P.A. and   Mansf\/ield E.L.,
{\it SIAM J. Appl. Math.}, 1994, V.55, 1693--1719.
\bibitem{gandarias:fu} Fushchych W. and Nikitin A.,
Symmetries of Maxwell's
Equations, Dordrecht, Reidel Publ. Comp., 1987.

\bibitem{gandarias:liu} Liu Y., {\it SIAM J. Math. Analysis,} 1995, V.26, 1527--1547.
\bibitem{gandarias:ml} Gandarias M.L., {\it J. Phys. A: Math. Gen.}, 1996, V.29,
607--633.
\bibitem{gandarias:ml1} Gandarias M.L., {\it J. Phys. A: Math. Gen.},
1996, V.29, 5919--5934.
\bibitem{gandarias:ml2} Gandarias M.L., {\it J. Phys. A: Math. Gen.}, 1997, V.30,
6081--6091.

\bibitem{gandarias:hill} Hill J.M., Dif\/ferential Equations and Group Methods,
CRC Press, Boca Raton, 1992.
\bibitem{gandarias:hi} Hill J.M., {\it J. Engng. Math.}, 1989, V.23, 141--155.

\bibitem{gandarias:hi1} Hill D.L. and  Hill J.M., {\it J. Engng. Math.},
 1990, V.24,  109--124.
\bibitem{gandarias:hi2} Hill J.M. and Hill D.L., {\it J. Engng. Math.},
1991, V.25, 287--299.
\bibitem{gandarias:ki1} King J.R., {\it J. Engng. Math.},
1988, V.22, 53--72.
\bibitem{gandarias:olv}  Olver P.J., Applications of Lie Groups to
Dif\/ferential Equations, Springer, New York, 1986.

\bibitem{gandarias:olv2} Olver P.J. and Rosenau P., {\it Phys. Lett. A.},
1986, V.114, 107--112.
\bibitem{gandarias:olv3} Olver P.J. and Rosenau P., {\it SIAM J. Appl. Math.},
 1987, V.47, 263--275.

\bibitem{gandarias:okra} Okrasinski,  {\it Extracta Mathematicae}, 1992, V.7,
93--95.
\bibitem{gandarias:ovsi} Ovsiannikov L.V.,
 Group Analysis of Dif\/ferential Equations,  Academic Press, New York, 1982.

\bibitem{gandarias:steph}  Stephani H.,  Dif\/ferential Equations:
Their Solution Using Symmetries, Cambridge U.P., Cambridge, 1989.

\end{thebibliography}
\end{document}